\documentclass[12pt,a4paper,oneside]{amsart}
\usepackage{amsfonts, amsmath, amssymb, amsthm, amscd}
\usepackage{hyperref}
\hypersetup{
  colorlinks   = true, 
  urlcolor     = blue, 
  linkcolor    = blue, 
  citecolor   = red 
}
\usepackage{anysize}
\marginsize{2cm}{2cm}{1.5cm}{1.5cm}

\usepackage[T2A]{fontenc}
\usepackage[cp1251]{inputenc}
\usepackage[english]{babel}
\usepackage{cmap}
\usepackage{enumerate}
\usepackage[pdftex]{graphicx}

\newtheorem{theorem}{Theorem}[section]
\newtheorem{lemma}[theorem]{Lemma}

\newtheorem{conjecture}[theorem]{Conjecture}
\newtheorem{corollary}[theorem]{Corollary}

\theoremstyle{definition}
\newtheorem{definition}[theorem]{Definition}

\theoremstyle{remark}
\newtheorem{remark}[theorem]{Remark}

\numberwithin{equation}{section}

\newcommand{\la}{\lambda}

\newcommand{\bra}{\langle}
\newcommand{\ket}{\rangle}
\newcommand{\R}{\mathbb{R}}

\def\EE{\mathcal{E}}%
\newcommand{\reff}[1]{(\ref{#1})}%

\newcommand{\ovl}[1]{\overline{#1}}

\def\Lin{\mathop{\rm Lin}}

\def\diag{\mathop{\rm diag}}

\!%
\def\vol{\mathop{\rm vol}}%

\newcommand{\crosp}{\Diamond}%
\newcommand{\cube}{\Box}%
\newcommand{\prooff}{\noindent {\bf Proof.}\\}%
\newcommand{\bbox}{\par\noindent\ensuremath{\Box}\par\noindent}%
\newenvironment{prf}
{\prooff}
{\bbox}

\begin{document}

\begin{center}
\textbf{\Large On the volume  of the John-L\"owner  ellipsoid}\\[0.3cm]

\emph{Grigory Ivanov\footnote{
DCG, FSB, Ecole Polytechnique F\'ed\'erale de Lausanne, Route Cantonale, 1015 Lausanne, Switzerland.
\\ \noindent
Department of Higher Mathematics, Moscow Institute of Physics and Technology,  Institutskii pereulok 9, Dolgoprudny, Moscow
region, 141700, Russia.
\\ \noindent
grimivanov@gmail.com
\\ \noindent 
Research partially supported by Swiss National Science Foundation grants 200020-165977 and 200021-162884
Supported by Russian Foundation for Basic Research, project 16-01-00259.
} 
}
\end{center}
\vspace{0.2cm}

{\bf Abstract.} We find an optimal upper bound on the volume of the John ellipsoid of a $k$-dimensional section of the $n$-dimensional cube, and an  optimal lower  bound on the volume of the L\"owner ellipsoid of a  projection of the $n$-dimensional cross-polytope onto a $k$-dimensional subspace. We use these results to  give a new proof of  Ball's upper bound on the volume of a $k$-dimensional section of the hypercube, and of Barthe's lower bound on the volume of a  projection of the $n$-dimensional cross-polytope onto a $k$-dimensional subspace. 
We settle equality cases in these inequalities. 
 Also, we describe all possible vectors in $\R^n,$ whose coordinates are the squared lengths of a projection of the standard basis in $\R^n$ onto a $k$-dimensional subspace. 
\bigskip


\section{Introduction}
In \cite{John}, Fritz John proved that each convex body in $\R^k$ contains a unique
ellipsoid of maximal volume.   John
characterized all convex bodies $K$ such that  the ellipsoid of maximal volume in $K$ is the Euclidean unit ball,  
$\EE_k.$ 

\begin{theorem}[F. John]
The Euclidean ball is the ellipsoid of maximal volume
contained in a convex body $K \in \R^k$ iff $\EE_k \subset K$ and, for some $n \geqslant k$,
there are Euclidean unit vectors $(u_i)_1^n $, on the boundary of $K$, and positive
numbers $(c_i)_1^n$ for which 

\begin{enumerate}
\item \label{John_cond_1}  $\sum\limits_1^n  c_i u_i = 0$;
\item \label{John_cond_2} $ \sum\limits_1^n c_i u_i \otimes u_i = I_k,$ the identity on $\R^k.$
\end{enumerate}
\end{theorem}

Keith Ball  in \cite{ball1992ellipsoids}  added the converse  part to this  theorem. 

\begin{theorem}[K. Ball]\label{th_ball}
 Let $(u_i)_1^n $ be a sequence of unit vectors in $R^k$ and $(c_i)_1^n$
be a sequence of positive numbers satisfying \reff{John_cond_1} and \reff{John_cond_2}.

Then the set $K = \{x \in \R^k | \bra x, u_i \ket  \leqslant 1,  i \in \ovl{1, n}\}$ contains a unique ellipsoid of
maximal volume, which is the Euclidean unit ball. 
\end{theorem}

The use of  vectors $(u_i)_1^n $ and positive numbers $(c_i)_1^n$ satisfying \reff{John_cond_1} and \reff{John_cond_2} appears to be extremely powerful in a range of problems in  convex analysis,
including (see \cite{ball2001convex} \nocite{zong2006cube}): 
tight  bounds on the volume ratio and on the outer volume ratio for centrally-symmetric  convex bodies, and optimal upper bounds on the volume of a $k$-dimensional section of the $n$-cube.

 In this short paper we study a simple alternative description of    the vectors $(u_i)_1^n \subset \R^k$ and positive numbers $(c_i)_1^n$ satisfying \reff{John_cond_2}.
\begin{definition}
We will say that some vectors $(v_i)_1^n \subset H$ {\it give  a unit decomposition} in  a $k$-dimensional vector space $H$ 
if 
\begin{equation}
\label{John_cond_3}
\sum\limits_1^n v_i \otimes v_i  = I_H,
\end{equation}
where $I_H$ is the identity operator in $H.$
\end{definition}

Clearly, non-zero vectors $(v_i)_1^n \subset \R^k$  give  a unit decomposition in $\R^k$   iff 
the vectors $\left(\frac{v_i}{|v_i|}\right)_1^n \subset \R^k$ and positive numbers $(|v_i|^2)_1^n$ satisfy \reff{John_cond_2}.

In Lemma \ref{lemma_realizable_vectors}, the set of  all possible vectors of positive reals $(c_1, \cdots, c_n),$ which we can get from condition \reff{John_cond_2}.
 That this result may be interesting for finding optimal bounds in different geometric inequalities, including the Brascamb--Lieb inequality. 

Using  a geometric approach, we will give in Theorem \ref{th_vol_ee} an optimal upper bound on the volume of the John ellipsoid of a $k$-dimensional section of the $n$-dimensional cube,  and derive an  optimal lower  bound on the volume of the L\"owner ellipsoid of a  projection of the $n$-dimensional cross-polytope onto a $k$-dimensional subspace. 

In Section \ref{section_last} we give a new proof   for  Ball's upper bound on the volume of a $k$-dimensional section of the hypercube  (see \cite{ball1989volumes}) and for  Barthe's lower bound on the volume of a  projection of the $n$-dimensional cross-polytope onto a $k$-dimensional subspace (see \cite{barthe1998reverse}). Moreover, we settle equality cases in these inequalities.

 We use $\bra p,x \ket$ to denote the \emph{value of a linear functional} $p$ \emph{at a vector}  $x.$  
For a  convex body  $K \subset \R^n$ we denote by $K^\circ$ and $\vol K$ the \emph{polar} body and the volume of $K$, respectively. 
We use $\crosp^n, \cube^n$ to denote the $n$-dimensional cross-polytope and cube, respectively. 
For a convex body $K \subset \R^n$ and a $k$-dimensional subspace $H_k$ of $\R^n$ we denote by
$K \cap H_k$ and $K | H_k$ the section of $K$ by $H_k$ and the orthogonal projection of $K$ onto $H_k,$ 
respectively. 

%
%

\section{Properties of a unit decomposition}

In this section, using a simple linear algebra observation  we introduce a  description  of  sets of  vectors $(v_i)_1^n \subset \R^k$ which give a unit decomposition in $\R^k$.  
Here we understand $\R^k \subset \R^n$ as the subspace of vectors, whose  last $n-k$ coordinates are zero.
For convenience, we will consider $(v_i)_1^n \subset \R^k \subset \R^n$ to be $k$-dimensional vectors. 



\begin{lemma}\label{lemma_equiv_cond}
The following assertions
 are equivalent:
\begin{enumerate}
\item \label{lemma_eq_cond_it1} vectors  $(v_i)_1^n \subset \R^k$  give a unit decomposition in $\R^k$; 
\item \label{lemma_eq_cond_it2} there exists an orthonormal basis $\{ f_1, \cdots, f_n \}$  in  $\R^n$ such that   $v_i$ is the orthogonal projection 
of $f_i$ onto $\R^k,$ for any $i \in \ovl{1,n};$ 
\item \label{lemma_eq_cond_Gram} 
$\Lin\{v_1, \cdots ,v_n\} = \R^k$ and the Gram matrix $\Gamma$ of vectors $\{v_1, \cdots ,v_n\} \subset \R^k$ is the matrix of a  projection operator from $\R^n$ onto the  linear hull of the rows the matrix 	$A = [v_1, \cdots, v_n].$
\item \label{lemma_eq_cond_it3}  the $k\times n$ matrix $A = [v_1, \cdots, v_n ]$ is a sub-matrix of an orthogonal matrix of order $n$.
\end{enumerate}	
\end{lemma}
\begin{prf}
1) \reff{lemma_eq_cond_it3} $\Rightarrow$ \reff{lemma_eq_cond_Gram}. \\ \noindent
Since $A$ is a sub-matrix of an orthogonal matrix, we have that $\operatorname{rk} A = k.$ Therefore, $\Lin\{v_1, \cdots ,v_n\} = \R^k.$

Let $\Gamma$ be the Gram matrix of  vectors $\{v_1, \cdots ,v_n\}$ and $P$ be the matrix 
 of a  projection operator from $\R^n$ onto the  linear hull of the rows of the  matrix 	$A = [v_1, \cdots, v_n].$ 
 Since the rows of $A$ form an orthonormal basis  of their linear hull $H_k$, 
 we can identify $Pe_i$ and $v_i$ in this basis of $H_k.$  Therefore, 
 $$
 	P_{ij} = \bra Pe_i, e_j\ket  = \bra P^2 e_i, e_j\ket  = \bra Pe_i, Pe_j\ket = \bra v_i, v_j \ket = \Gamma_{ij}.
 $$
\\

\noindent
2) \reff{lemma_eq_cond_Gram} $\Rightarrow$ \reff{lemma_eq_cond_it2}. \\ \noindent
Let the Gram matrix $\Gamma$ of the vectors $\{v_1, \cdots ,v_n\}$ be the matrix of a  projection operator onto $H_k$.
By the last identity, we have that $\bra \Gamma e_i,  \Gamma e_j\ket = \bra v_i, v_j \ket.$
But if two sets $S_1$ and $S_2$ of vectors have the same Gram matrix, then there exists an orthogonal transformation of the space that maps 
vectors of  $S_1$ to $S_2$. Indeed, each step in the  Gram-Schmidt process  for both systems are identical, that means that any orthogonal transformation which maps  the Gram-Schmidt orthonormalization of $S_1$ to the Gram-Schmidt orthonormalization of $S_2$ maps $S_1$ to $S_2.$  
Therefore, with a proper choice of the orthonormal basis in $\R^n,$ we can identify vectors $v_i$ and $Pe_i,$ for $i \in \ovl{1,n},$ ans subspaces $H_k$ and $\R^k = \Lin\{v_1, \cdots ,v_n\}.$      

\noindent
3)\reff{lemma_eq_cond_it2}  $ \Rightarrow$ \reff{lemma_eq_cond_it1}\\ \noindent
Let $P$ be the projection from  $\R^n$ onto $\R^k.$ Let  $v_i = Pf_i.$ 
For an arbitrary vector $x \in \R^k$ we have $Px = x$ and therefore  
$$ \left(\sum\limits_1^n v_i \otimes v_i\right) x = \sum\limits_1^n \bra v_i, x \ket v_i  = 
\sum\limits_1^n \bra Pf_i, x \ket v_i = \sum\limits_1^n \bra f_i, P x \ket v_i = 
P(\sum\limits_1^n \bra f_i, x \ket f_i) = Px = x.$$ \\

\noindent
4)\reff{lemma_eq_cond_it1}  $\Rightarrow$ \reff{lemma_eq_cond_it3}\\ \noindent
For $i \in \ovl{1,k}$ we have $e_j = \sum\limits_1^n \bra v_i, e_j \ket v_i.$
Therefore,
$$
	1 = |e_j|^2 = \sum\limits_1^n \bra v_i, e_j \ket^2 = \sum\limits_1^n v_i[j]^2,
$$ 
a2nd 
$$
	0 = \bra e_j, e_k \ket =  \sum\limits_1^n \bra v_i, e_j \ket \bra v_i, e_k \ket =
	\sum\limits_1^n  v_i[j] v_i[k],
$$
where $v_i[j]$	 is the $j$'th coordinate of the vector $v_i$ in the given basis.

That is, the rows of the $k \times n$ matrix $[v_1, \cdots, v_n ]$ form an orthonormal system of $k$ vectors in $\R^n$. 
\end{prf}


As a direct consequence of Lemma \ref{lemma_equiv_cond} we get
\begin{corollary}	
Let $(u_i)_1^n $ be a sequence of unit vectors in $R^k$ and $(c_i)_1^n$
be a sequence of positive numbers satisfying  \reff{John_cond_2}.
Then the set $K = \{x \in \R^k | |\bra x, u_i \ket|  \leqslant 1,  i \in \ovl{1, n}\}$ is an affine image of a $k$-dimensional section of $\cube^n.$
\end{corollary}

\begin{definition}
We will say that a  vector $C = (c_1, \cdots, c_n)$ is {\it realizable} in $\R^k$ if there exist   vectors  
$(v_i)_1^n$ which give a unit decomposition in $\R^k$ such that $c_i = |v_i|^2,$ $i \in \ovl{1,n}.$
\end{definition} 

Now we are going to describe all possible realizable vectors in $\R^k$.  
For this purpose, we need to use the following standard notation.

\begin{definition}
Let $a$ and $b$ be non-negative vectors in $\R^n$.
The vector $a$   {\it majorizes} the vector $b,$ which we denote by $a \succ b,$ 
if the sum of the $k$ largest entries of $a$ is at least   the sum of the $k$ largest entries of $b,$ 
for every $k \in \ovl{1,n},$ and the sums of all entries of $a$ and $b$ are equal.
\end{definition}

\begin{lemma} \label{lemma_realizable_vectors}
	A vector $(c_1, \cdots, c_n)$ is realizable iff 
	$$(\underbrace{1, \cdots, 1}_k , \underbrace{ 0, \cdots, 0}_{n-k}) \succ (c_1, \cdots, c_n). $$
\end{lemma}
\begin{prf}
Let $(c_1, \cdots, c_n)$ be a realizable vector. By definition and by Lemma \ref{lemma_equiv_cond}, there are vectors $(v_i)_1^n \subset \R^k$ that give a unit decomposition in $\R^k$ such that   the diagonal entries of their Gram matrix $\Gamma$ are   $(c_i)_1^n,$ and 
$\Gamma$ is the matrix of a projection operator from $\R^n$ onto some $k$-dimensional subspace $H_k.$   

 So the vector $(c_1, \cdots, c_n)$ is realizable iff  there exists a projection operator from $\R^n$ onto some $k$-dimensional subspace   with $(c_1, \cdots, c_n)$ on the main diagonal. Applying Horn's theorem \cite{horn1954doubly}, 
 which asserts that a vector $(c_1, \cdots, c_n)$ can be the main diagonal of a Hermitian matrix with a vector of eigenvalues 
  $(\la_1, \cdots, \la_n)$ iff  $\la \succ c,$  to the  vector 
  $ (\underbrace{1, \cdots, 1}_k , \underbrace{ 0, \cdots, 0}_{n-k})$ we complete the proof.
\end{prf}

\section{Estimation for the volume of the L\"owner-John ellipsoid}

Before stating the next result, we recall that the John ellipsoid of a convex body $K$ is the ellipsoid of maximal volume contained in $K$, and the L\"owner ellipsoid of a convex body $K$ is the ellipsoid of minimal volume containing $K.$ We use  $\mathcal{E}_{H_k}$ and $\EE^{H_k}$ to denote the L\"owner ellipsoid of  $\crosp^n | H_k$ and the John ellipsoid of  $\cube^n \cap H_k,$ respectively.

\begin{lemma}\label{lemma_covering}
	Suppose  vectors $(v_i)_1^n \subset \R^k$  give a unit decomposition in $\R^k$.  Then for any ellipsoid $\EE$ with the center in the origin that covers all vectors $(v_i)_1^n,$  we have 
\begin{equation} \label{lemma_covering_unit_dec}
\vol \EE \geqslant \left(\frac{k}{n}\right)^{\frac{k}{2}} \vol \EE_k,
\end{equation}	
where $\EE_k$ is the unit ball in $\R^k.$
The bound is tight.
The inequality becomes an equality  iff $|v_i|^2 =\frac{k}{n},$ for all $i \in \ovl{1,n}.$
\end{lemma}
\begin{prf}
	For a  positive-definite operator $A$ on $\R^k,$ we know that the volume of the ellipsoid
	$\{x \in \R^k| \bra Ax, x \ket \leqslant 1\}$ is $\frac{\vol \EE_k}{\sqrt{\det A}}.$
	To prove our lemma, it is enough to show that for any positive-definite  operator $A$ such that 
	$\bra Av_i, v_i \ket \leqslant 1$ for $i \in \ovl{1,n}$ we have  $\det A \leqslant \left(\frac{n}{k}\right)^k.$
	
	Fix a  positive-definite operator $A$ on $\R^k$ such that  $\bra Av_i, v_i \ket \leqslant 1$ for $i \in \ovl{1,n}.$
	We can choose an orthonormal basis in $\R^k$ such that $A = \diag \{\la_1, \cdots, \la_k\}$ in this basis. 
	Let $v'_i$ be the coordinate vector of $v_i$  in the new basis for each  $i\in \ovl{1,n}$.  
	We can rewrite the inequality $\bra Av_i, v_i \ket \leqslant 1$ in the following form:
	\begin{equation}\label{lemma_cover_unit_dec_2}
		\sum\limits_1^k \la_j v'_i[j]^2 \leqslant 1,
	\end{equation}
where $v'_i[j]$	 is the $j$'th coordinate of the vector $v'_i$ in the given basis. 
	
		Summing up the  inequalities (\ref{lemma_cover_unit_dec_2}) for all $i \in \ovl{1,n}$ and using the observation that 
	$\sum\limits_1^n  v'_i[j]^2 = 1 $  (see Lemma  \ref{lemma_equiv_cond}), we get
	$$
		\sum\limits_1^k  \la_i \leqslant n.
	$$
	Applying the  AM-GM inequality, we get
	$$
		\det A = \prod\limits_1^k \la_i \leqslant
		 \left( \frac{\sum\limits_1^k  \la_i}{k} \right)^k \leqslant \left( \frac{n}{k} \right)^k.
	$$ 
	
	According to Lemma \ref{lemma_realizable_vectors}, the vector $\left(\frac{k}{n}, \cdots,  \frac{k}{n} \right)$ 
	is realizable. It is clear that in this case inequality \reff{lemma_covering_unit_dec} becomes an equality.
	Moreover, by properties of the AM-GM inequality, we have an equality in \reff{lemma_covering_unit_dec} iff $\la_i = \frac{n}{k},$ for all $i \in \ovl{1,n}.$ 
	By the inequality \reff{lemma_cover_unit_dec_2}, there is an equality  in \reff{lemma_covering_unit_dec}  iff $|v_i|^2 =\frac{k}{n},$ for all $i \in \ovl{1,n}.$ 
\end{prf}

Fix a $k$-dimensional subspace $H_k$ in $\R^n.$ Let $P : \R^n \to H_k$ be the projection onto $H_k.$
Since $\crosp^n | H_k$ is the absolute convex hull of the vectors 
$v_i = Pe_i$ for $i \in \ovl{1,n}$ that give us a unit decomposition in $H_k$, Lemma \ref{lemma_covering}  implies 
that the volume of the L\"owner ellipsoid for $\crosp^n | H_k$  is at least 
$
\left(\frac{k}{n}\right)^{\frac{k}{2}} \vol \EE_k.
$ 

To settle the reverse case, we need to recall the following simple duality arguments.

For a given $k$-dimensional subspace $H_k$ in  $\R^n,$ we can consider 
 the space $H_k \subset (\R^n)^* = \R^n$ itself to be the dual space for $H_k$.  
 Indeed, $H_k$ is a $k$-dimensional space consisting of all linear functionals on $H_k$ with the proper linear structure, and the restriction of   the Euclidean norm in $\R^n$  onto $H_k$ generates   the operator norm on $H_k.$

For the sake of completeness we give a proof of the following well-known 
fact.
\begin{lemma}\label{lemma_dual}
Let $H_k$ be a $k$-dimensional subspace of $\R^n.$ Assume the dual space $H_k^*$ for $H_k$ is $H_k$ itself. 
For a convex body $K \in \R^n$ containing the origin in the interior, we have 
	$$
		(K \cap H_k)^\circ = K^\circ | H_k,
	$$
where we understand $K \cap H_k$	as a subset of $H_k,$ and its polar set as a subset of $H_k^* = H_k.$  
\end{lemma}
\begin{prf}
We  use  $H_k^\perp$ to denote the orthogonal complement of $H_k$ in $\R^n.$

1) Let us show that $(K \cap H_k)^\circ \supset K^\circ | H_k.$ 
Fix a functional $p \in K^\circ  | H_k.$ Since $p$ belongs to the projection of $K^\circ$, there is  a functional 
$p^\perp \subset H_k^\perp$ such that $p + p^\perp \in K^\circ.$ By definition of the polar body, we have 
$ \bra p + p^\perp, x \ket \leqslant 1$ for any $x \in K.$  In particular, for any $x \in K \cap H_k,$ we have
$$
	1 \geqslant \bra p + p^\perp, x \ket = \bra p, x \ket  +  \bra p^\perp, x \ket = \bra p, x\ket.
$$   
This means that  $p \in (K \cap H_k)^\circ.$

2)  Let us show that $(K \cap H_k)^\circ \subset K^\circ | H_k.$ 
Suppose for a contradiction that  there is a functional 
$p  \in \left(K \cap H_k \right)^\circ$ such that $p \notin K^\circ | H_k.$ 
By the hyperplane separation theorem, there exists a vector $y \in H_k$ such that
\begin{equation}\label{lemma_dual_1}
\bra p, y \ket > 1 
\end{equation}	
and
\begin{equation}\label{lemma_dual_2}
\bra q, y \ket \leqslant 1  
	\quad \mbox{for all}\quad  q \in   
	 K^\circ  | H_k.
\end{equation}
Clearly, $\bra y , q^\perp\ket = 0 $ for any $q^\perp \in H_k^\perp.$
Combining this and the inequality \reff{lemma_dual_2}, we get 
$$
	\bra q, y \ket \leqslant 1 
$$
for all $q \in K^\circ.$ 
By the definition of the polar set, we obtain 
$y \in \left(K^\circ \right)^\circ = K.$ 
So $y \in K$ and $y \in H_k,$ therefore $y \in K \cap H_k.$
This contradicts  the inequality \reff{lemma_dual_1}. 
\end{prf}

For a given  convex centrally-symmetric body $K$ with the center at the origin, by symmetry and duality arguments, we have that the polar ellipsoid of the John ellipsoid of $K$ is the L\"owner ellipsoid of $K^\circ.$ 
 
Summarizing the arguments of section 3, we obtain.

\begin{theorem}	\label{th_vol_ee}
For any $1 \leqslant k \leqslant n$ we have
$$
	\frac{\vol \EE_{H_k}}{\vol \EE_k} \geqslant \left(\frac{k}{n}\right)^{\frac{k}{2}} \mbox{  and  }\quad  
	\frac{\vol \EE^{H_k}}{\vol \EE_k} \leqslant \left(\frac{n}{k}\right)^{\frac{k}{2}}. 
$$
The bounds are sharp. That is, there exists a subspace $H_k$ such that the two inequalities are simultaneously
 hold as equalities.  
\end{theorem}
\section{Bounds on  the volume of a section of $\cube^n$ and a projection of $\crosp^n$}\label{section_last}
K. Ball, in his fundamental paper \cite{ball1989volumes}, proved the following inequality
\begin{equation}\label{ineq_Ball_1}
\frac{\vol (\cube^n \cap H_k)}{ \vol \cube^k} \leqslant \left(\frac{n}{k}\right)^{\frac{k}{2}}.
\end{equation}
F. Barthe in \cite{barthe1998reverse} proved the dual inequality
\begin{equation}\label{ineq_Barthe_1}
\frac{\vol (\crosp^n | H_k)}{ \vol \crosp^k} \geqslant \left(\frac{k}{n}\right)^{\frac{k}{2}}.
\end{equation}
One can see that both inequalities become equalities when $k | n$ and $H_k$ is determined by the system of linear equations 
\begin{equation}\label{equality_case}
x_{\frac{n}{k} j + i_1} = x_{\frac{n}{k} j + i_2}, \mbox{  where}\quad j \in \ovl{0, k-1} \mbox{  and} \quad 
1 \leqslant i_1, i_2  \leqslant  \frac{n}{k}.
\end{equation}

 Using Theorem \ref{th_vol_ee}, we are going to  give another proof of 
the inequalities \reff{ineq_Ball_1} and \reff{ineq_Barthe_1}, and settle the equality case. 
\begin{theorem} \label{corollary}
For any $k$-dimensional subspace of $\R^n,$ we have
	$$
	\frac{\vol (\cube^n \cap H_k)}{ \vol \cube^k} \leqslant \left(\frac{n}{k}\right)^{\frac{k}{2}} 
	\mbox{  and  }\quad   
	\frac{\vol (\crosp^n |  H_k)}{ \vol \crosp^k} \geqslant \left(\frac{k}{n}\right)^{\frac{k}{2}}.
	$$
The bounds are optimal iff $k|n$.
\end{theorem}	
\begin{prf}
Using the Brascamb--Lieb inequality, K. Ball \cite{ball1989volumes} proved that among all $k$-dimensional convex centrally-symmetric bodies,
the $k$-cube has the greatest volume ratio 
(i.e., $\left( \frac{\vol{K}}{\vol \EE} \right)^{\frac{1}{n}}$, where $\EE$ is the John ellipsoid of $K$). This means that
\begin{equation}\label{ineq_vol_ratio}
	\frac{\vol (\cube^n \cap H_k)}{ \vol \cube^k} \leqslant \frac{\vol \EE^{H_k}}{ \vol \EE_k}.
\end{equation}

The dual case for the outer volume ratio (i.e. $\left( \frac{\vol \EE}{\vol K} \right)^{\frac{1}{n}}$, where $\EE$ is the L\"owner ellipsoid of $K$) was resolved using Barthe's reverse Brascamb--Lieb inequality \cite{ball2001convex}.
It was shown that $\crosp^k$ has the biggest outer volume ratio among all $k$-dimensional convex centrally-symmetric bodies. Therefore 
\begin{equation}\label{ineq_invol_ratio}
	\frac{\vol (\crosp^n | H_k)}{ \vol \crosp^k} \geqslant \frac{\vol \EE_{H_k}}{ \vol \EE_k}.
\end{equation}

Combining \reff{ineq_vol_ratio} and \reff{ineq_invol_ratio} with the inequalities from the assertion of Theorem \ref{th_vol_ee}, we obtain \reff{ineq_Ball_1} and \reff{ineq_Barthe_1}.

We now prove that the bounds are optimal only if $k|n.$

In \cite{barthe1998reverse} Proposition 10, Barthe proved that whenever the volume ratio for a convex centrally-symmetric body $K \subset \R^k$ equals 
the volume ratio for $\cube^k,$ then $K$ is an affine $k$-dimensional cube (or parallelotope).  Also, he proved that 
if a  centrally-symmetric convex  body $K \subset \R^k$ has the extremal inner volume ratio, then $K$ is an  affine cross-polytope. These arguments imply that 
$\cube^n \cap H_k$ is an affine cube   and  $\crosp^n | H_k$ is an affine cross-polytope in the equality cases for the inequalities \reff{ineq_Ball_1} and \reff{ineq_Barthe_1}, respectively.

 Using the fact that $\cube^k$ is the polar set of $\crosp^k$ and employing Lemma \ref{lemma_dual}, we obtain that 
 for any given subspace $H_k$  equality holds in  \reff{ineq_Ball_1} if and only if  equality holds in  \reff{ineq_Barthe_1}.
 Hence,  it is enough to settle  equality  only for the inequality \reff{ineq_Barthe_1}.
 
 Suppose for a given $H_k$ we have  equality in \reff{ineq_Barthe_1}.   Then $\crosp^n | H_k$ is an affine $k$-dimensional cross-polytope. Let $P$ be the projection from  $\R^n$ onto $H_k,$ and   $v_i = Pe_i.$ 
It is easy to see that   each vertex of $\crosp^n | H_k$  is identical to at least one of the vectors $v_i,$ where $i \in \ovl{1,n}.$
  The proof of Lemma \ref{lemma_covering} yields that all  lengths $|v_i|,$ for $i \in \ovl{1,n},$ are the same. From this and the triangle inequality,
we conclude that all  vectors $v_i,$ $i \in \ovl{1,n},$ are vertices of the affine cross-polytope  $\crosp^n | H_k.$
So, each vertex of $\crosp^n|H_k$ is identical to some $v_i$, and conversely, each $v_i$ is a vertex of $\crosp^n|H_k$.

Denote by $\ell_i,$ $i \in \ovl{1,k},$ lines in $H_k$  that pass through vertices of the affine cross-polytope   
$\crosp^n | H_k.$  We showed that for any $i \in \ovl{1,n}$ there exists $j \in \ovl{1,k}$ such that $v_i \in \ell_j.$ 
Hence, there exist vectors $d_j \in \ell_j,$ $j \in \ovl{1,k},$ such that  
$I_{H_k} = \sum\limits_{1}^{n}v_i \otimes v_i =  \sum\limits_{1}^{k}d_j \otimes d_j.$  
This means that the vectors $d_j,$ $j \in \ovl{1,k},$ give us a   unit decomposition in $H_k.$
By the assertion \ref{lemma_eq_cond_it3} of Lemma \ref{lemma_equiv_cond}, we have that the vectors 
$d_j,$ $j \in \ovl{1,k},$  form an orthonormal basis in $H_k.$  
Therefore, all $k$ sums  $\sum\limits_{v_i \in \ell_j} |v_i|^2,$ $j \in \ovl{1,k},$   equals 1.
As mentioned above, all lengths $|v_i|,$ for $i \in \ovl{1,n},$ are the same. 
 Consequently, the same number of vectors $v_i,$  $i \in \ovl{1,n},$  lies on each line $\ell_j,$ $j \in \ovl{1,k}.$  That is, $k | n.$   
\end{prf}
\begin{remark}
	Up to coordinate permutation and up to change of the sign of coordinates, Theorem \ref{corollary} implies that  equality  in \reff{ineq_Ball_1} and \reff{ineq_Barthe_1}  is attained when  $H_k$ is determined by  \reff{equality_case}.
\end{remark}

We should note  that  Ball's and Barthe's proofs of the inequalities \reff{ineq_Ball_1} and \reff{ineq_Barthe_1}
used the same arguments as the proofs of \reff{ineq_vol_ratio} and \reff{ineq_invol_ratio}. 
However, we believe that it  may be of interest how our result reveals the connection between  Theorem 
\ref{corollary} and the volume of the L\"owner and the John ellipsoid.

We conjecture the following.
\begin{conjecture}
\begin{equation}\label{conjecture}
	\frac{\vol (\crosp^n |  H_k)}{ \vol \crosp^k} \geqslant {2}^{\frac{k-n}{2}}.
\end{equation}

The bound is optimal when $2k \geqslant n.$
\end{conjecture}
This is the dual statement for  another Ball's upper bound on the volume of a $k$-dimensional section of $\cube^n:$
\begin{equation}\label{ineq_Ball_2}
	\frac{\vol (\cube^n \cap H_k)}{ \vol \cube^k} \leqslant {2}^{\frac{n-k}{2}}.
\end{equation}

We note that  inequalities \reff{ineq_Barthe_1} and  \reff{conjecture}  follow from the well-known Mahler conjecture and  inequalities \reff{ineq_Ball_1} and \reff{ineq_Ball_2}, respectively.

{\bf Acknowledgements.}  
I am grateful for Marton Naszodi for useful remarks and help with the text.
I am especially grateful for Roma Karasev for fruitful conversations concerned the topic of this paper. 

\nocite{brazitikos2014geometry}
\bibliographystyle{abbrv}
\bibliography{/home/grigory/Dropbox/work_current/uvolit} 
\end{document}